\documentclass[preprint]{aastex}
\usepackage{amssymb,amsmath,graphicx,natbib}

\newcommand{\ds}{\displaystyle}
\newcommand{\bc}{\begin{equation*}\begin{array}{l}}
\newcommand{\ec}{\end{array}\end{equation*}}
\newcommand{\ason}{\renewcommand{\arraystretch}{1.5}}
\newcommand{\asoff}{\renewcommand{\arraystretch}{1.0}}

\shorttitle{Multistep Methods with Improved Phase Lag Characteristics }
\shortauthors{D. S. Vlachos, Z. A. Anastassi and T. E. Simos}
\begin{document}

\title{A New Family of Multistep Methods with Improved Phase Lag Characteristics for the Integration of Orbital Problems}
\author{D. S. Vlachos\altaffilmark{1}, Z. A. Anastassi\altaffilmark{2}, and T. E. Simos\altaffilmark{3}}
\affil{Laboratory of Computational Sciences,\\
Department of Computer Science and Technology,\\
Faculty of Sciences and Technology,\\
University of Peloponnese, GR-22 100, Tripolis, Greece}

\altaffiltext{1}{e-mail: dvlachos@uop.gr}
\altaffiltext{2}{e-mail: zackanas@uop.gr}
\altaffiltext{3}{Highly Cited Researcher, Active Member of the European Academy of Sciences and Arts, Address: Dr. T.E. Simos, 26 Menelaou Street, Amfithea - Paleon Faliron, GR-175 64 Athens, GREECE, Tel: 0030 210 94 20 091, e-mail: tsimos.conf@gmail.com, tsimos@mail.ariadne-t.gr}

\begin{abstract}
In this work we introduce a new family of ten-step linear multistep methods for the integration of orbital problems. The new methods are constructed by adopting a new methodology which improves the phase lag characteristics by vanishing both the phase lag function and its first derivatives at a specific frequency. The efficiency of the new family of methods is proved via error analysis and numerical applications.

\end{abstract}


\section{Introduction}
The numerical integration of systems of ordinary differential equations with oscillatory solutions has been the subject of research during the past decades. This type of ODEs is often met in real problems, like the N-body problem. For highly oscillatory problems standard non-specialized methods can require a huge number of steps to track the oscillations. One way to obtain a more efficient integration process is to construct numerical methods with an increased algebraic order, although the implementation of high algebraic order meets several difficulties \cite{quinlan_arxiv_astro_ph_9901136}.

On the other hand, there are some special techniques for optimizing numerical methods. Trigonometrical fitting and phase-fitting are some of them, producing methods with variable coefficients, which depend on $v = \omega h$, where $\omega$ is the dominant frequency of the problem and $h$ is the step length of integration. More precisely, the coefficients of a general linear method are found from the requirement that it integrates exactly powers up to degree $p+1$. For problems having oscillatory solutions, more efficient methods are obtained when they are exact for every linear combination of functions from the reference set
\begin{equation}
\{1, x, \ldots , x^K , e^{\pm \mu x},\ldots , x^P e^{\pm \mu x}\}\label{equ_exp_fit}
\end{equation}
This technique is known as exponential (or trigonometric if $\mu=i\omega$) fitting and has a long history \cite{gautschi_NM_3_381_61}, \cite{lyche_NM_19_65_72}. The set (\ref{equ_exp_fit}) is characterized by two integer parameters, $K$ and $P$ . The set in which there is no classical component is identified by $K =-1$ while the set in which there is no exponential fitting component (the classical case) is identified by $P =-1$. Parameter $P$ will be called the level of tuning. An important property of exponential fitted algorithms is that they tend to the classical ones when the involved frequencies tend to zero, a fact which allows to say that exponential fitting represents a natural extension of the classical polynomial fitting. The examination of the convergence of exponentially fitted multistep methods is included in Lyche's theory \cite{lyche_NM_19_65_72}. There is a large number of significant methods presented with high practical importance that have been presented in the bibliography (see for example \cite{simos_Book_CMATV1_RCS_00}, \cite{chawla_JCAM_15_329_86}, \cite{raptis_CPC_14_1_78}, \cite{anastassi_IJMPC_15_1_04}, \cite{anastassi_MCM_42_877_05}, \cite{anastassi_JMC_41_79_07}, \cite{anastassi_NA_10_301_05}, \cite{lambert_JIMA_18_189_76}, \cite{cash_JNAIAM_1_81_06}, \cite{iavernaro_JNAIAM_1_91_06}, \cite{mazzia_JNAIAM_1_131_06}, \cite{berghe_JNAIAM_1_241_06}, \cite{psihoyios_CL_2_51_06}, \cite{simos_CL_1_37_05}, \cite{simos_CL_1_45_07}. The general theory is presented in detail in \cite{ixaru_Book_EF_KAP_04}.

Considering the accuracy of a method when solving oscillatory problems, it is more appropriate to work with the phase-lag, rather than the principal local truncation error. We mention the pioneering paper of Brusa and Nigro \cite{brusa_IJNME_15_685_80}, in which the phase-lag property was introduced. This is actually another type of a truncation error, i.e. the angle between the analytical solution and the numerical solution. On the other hand, exponential fitting is accurate only when a good estimate of the dominant frequency of the solution is known in advance. This means that in practice, if a small change in the dominant frequency is introduced, the efficiency of the method can be dramatically altered. It is well known, that for equations similar to the harmonic oscillator, the most efficient exponential fitted methods are those with the highest tuning level.

In this paper we present a new family of methods based on the 10-step linear multistep method of Quinlan and Tremain \cite{quinlan_AJ_100_1694_90}. The new methods are constructed by vanishing the phase-lag function and its first derivatives at a predefined frequency. Error analysis and numerical experiments show that the new methods exhibit improved characteristics concerning the long term behavior of the solution in the 5-outer problem. The paper is organized as follows: In section 2, the general theory of the new methodology is presented. In section 3, the new methods are described in detail. In section 5 the stability properties of the new methods are investigated. Section 5 presents the results from the numerical experiments and finally, conclusions are drawn in section 6.

\section{Phase-lag analysis of symmetric multistep methods}
Consider the differential equations
\begin{equation}
\frac{d^2y(t)}{dt^2}=f(t,y),\;y(t_0)=y_0,\;y'(t_0)=y'_0
\label{equ_ref_diff}
\end{equation}
and the linear multistep methods
\begin{equation}
\sum_{j=0}^{J}a_jy_{n+j}=h^2\sum_{j=0}^{J}b_jf_{n+j}
\label{equ_ref_meth}
\end{equation}
where $y_{n+j}=y(t_0+(n+j)h)$, $f_{n+j}=f(t_0+(n+j)h,y(t_0+(n+j)h))$ and $h$ is the step size of the method. We associate the following functional to method (\ref{equ_ref_meth}):
\begin{equation}
L(h,a,b,y(t))=\sum_{j=0}^Ja_jy(t+j\cdot h)-h^2\sum_{j=0}^Jb_jy''(t+j\cdot h)
\end{equation}
where $a,b$ are the vectors of coefficients $a_j$ and $b_j$ respectively, and $y(t)$ is an arbitrary function. The algebraic order of the method (\ref{equ_ref_meth}) is $p$, if
\begin{equation}
L(h,a,b,y(t))=C_{p+2}h^{p+2}y^{(p+2)}(t)+O(h^{p+3})
\label{equ_ref_error}
\end{equation}
The coefficients $C_q$ are given
\begin{eqnarray}
C_0=&\sum_{j=0}^J a_j \nonumber \\
C_1=&\sum_{j=0}^J j\cdot a_j \nonumber \\
C_q=&\frac{1}{q!}\sum_{j=0}^Jj^q\cdot a_j -\frac{1}{(q-2)!}\sum_{j=0}^Jj^{q-2}b_j
\end{eqnarray}
The principal local truncation error (PLTE) is the leading term of (\ref{equ_ref_error})
\begin{equation}
 plte=C_{p+2}h^{p+2}y^{(p+2)}(t)
\end{equation}
The following assumptions will be considered in the rest of the paper:
\begin{enumerate}
\item $a_J=1$ since we can always divide the coefficients of (\ref{equ_ref_meth}) with $a_J$.
\item $|a_0|+|b_0|\neq 0$ since otherwise we can assume that $J=J-1$.
\item $\sum_{j=0}^J |b_j| \neq 0$ since otherwise the solution of (\ref{equ_ref_meth}) would be independent of (\ref{equ_ref_diff}).
\item The method (\ref{equ_ref_meth}) is at least of order one.
\item The method (\ref{equ_ref_meth}) is zero stable, which means that the roots of the polynomial
\begin{equation}
p(z)=\sum_{j=0}^Ja_jz^j
\end{equation}
all lie in the unit disc, and those that lie on the unit circle have multiplicity one.
\item The method (\ref{equ_ref_meth}) is symmetric, which means that
\begin{equation}
a_j=a_{J-j},\;b_j=b_{J-j},\;j=0(1)J
\end{equation}
It is easily proved then that both the order of the method and the step number $J$ are even numbers \cite{lambert_JIMA_18_189_76}.
\end{enumerate}
Consider now the test problem
\begin{equation}
y''(t)=-\omega ^2 y(t)
\label{equ_ref_equ}
\end{equation}
where $\omega$ is a constant. The numerical solution of (\ref{equ_ref_equ}) by applying method (\ref{equ_ref_meth}) is described by the difference equation
\begin{equation}
\sum_{j=1}^{J/2} A_j(s^2)(y_{n+j}+y_{n-j})+A_0(s^2)y_n=0
\end{equation}
with
\begin{equation}
A_j(s^2)=a_{\frac{J}{2}-j}+s^2\cdot b_{\frac{J}{2}-j}
\end{equation}
and $s=\omega h$. The characteristic equation is then given
\begin{equation}
\sum_{j=1}^{J/2} A_j(s^2)(z^j+z^{-j})+A_0(s^2)=0
\label{equ_char_equ}
\end{equation}
and the interval of periodicity $(0,s_0^2)$ is then defined such that for $s\in (0,s_0)$ the roots of (\ref{equ_char_equ}) are of the form
\begin{equation}
z_1=e^{i\lambda (s)},\;z_2=e^{-i\lambda (s)},\;|z_j|\leq 1,\;3\leq j\leq J
\end{equation}
where $\lambda (s)$ is a real function of $s$. The phase-lag $PL$ of the method (\ref{equ_ref_meth}) is then defined
\begin{equation}
PL=s-\lambda (s)
\end{equation}
and is of order $q$ if
\begin{equation}
PL=c\cdot s^{q+2}+O(s^{q+4})
\end{equation}
In general, the coefficients of the method (\ref{equ_ref_meth}) depend on some parameter $v$, thus the coefficients $A_j$ are functions of both $s^2$ and $v$. The following theorem has been proved by Simos and Williams \cite{simos_CC_23_513_99}:
For the symmetric method (\ref{equ_ref_equ}) the phase-lag is given by
\begin{equation}
PL(s,v)=\frac{2\sum_{j=1}^{J/2}A_j(s^2,v)\cdot cos(j\cdot s) +A_0(s^2,v)}{2\sum_{j=1}^{J/2}j^2A_j(s^2,v)}
\end{equation}
We are now in position to describe the new methodology. In order to efficiently integrate oscillatory problems, a good technique is to calculate the coefficients of the numerical method by forcing the phase lag to be zero at a specific frequency. But, since the appropriate frequency is problem dependent and in general is not always known, we may assume that we have an error in the frequency estimation. It would be of great importance to force the phase-lag to be less sensitive to this error. Thus, beyond the vanishing of the phase-lag, we also force its first derivatives to be zero.

\section{Construction of the new methods}
\subsection{Base Method}
The family of new methods is based on the 10-step linear multistep method of Quinlan and Tremain \cite{quinlan_AJ_100_1694_90} which is of the form (\ref{equ_ref_meth}) with coefficients
\begin{equation}
\begin{array}{cccccc}
a_0=1 & a_1=-1 & a_2=1 & a_3=-1 & a_4=1 & a_5=-2 \\
b_0=0 & b_1=\frac{399187}{241920} & b_2=-\frac{17327}{8640} & b_3=\frac{597859}{60480} & b_4=-\frac{704183}{60480} & b_5=\frac{465133}{24192}
\end{array}
\label{equ_base_meth}
\end{equation}
The principal term of the local truncation error (PLTE) of the method is given
\begin{equation}
 plte=\frac {52559}{912384}\,y^{(12)} {h}^{12}
\end{equation}
\subsection{Method PF-D0: Phase fitted}
The first method of the family (PF-D0) is constructed by forcing the phase-lag function to be zero at the frequency $V=\omega *h$. Coefficients $a$ are left the same, while coefficients $b$ become:

\begin{equation*}
\begin{array}{rlrlrl}
\ds b_0=&0 , & b_1=&\frac{1}{8064}\frac{b_{1,num}^0}{D_0} , & b_2=&\frac{1}{4032}\frac{b_{2,num}^0}{D_0} , \\
\ds b_3=&\frac{1}{2016}\frac{b_{3,num}^0}{D_0} , & b_4=&\frac{1}{4032}\frac{b_{4,num}^0}{D_0} , & b_5=&\frac{1}{4032}\frac{b_{5,num}^0}{D_0}
\end{array}
\end{equation*}

where
\begin{equation}
D_0={v}^{2} ( ( \cos ( v ) ) ^{4}+6\, (
\cos ( v ) ) ^{2}+1-4\,\cos ( v ) -4\,
 ( \cos ( v ) ) ^{3} )
\nonumber
\end{equation}
and
\begin{equation*}
\begin{array}{l}
b_{1,num}^0=-16128\,(\cos(v))^{3}+45139\,(\cos(v))^{3}{v}^{2}+6048\,(\cos(v))^{2}-73215\,(\cos(v))^{2}{v}^{2}\\
+3024\,\cos(v)+47553\,\cos(v){v}^{2}-11917\,{v}^{2}+16128\,(\cos(v))^{5}-8064\,(\cos(v))^{4}-1008
\\
b_{2,num}^0=-32256\,(\cos(v))^{4}+45139\,(\cos(v))^{4}{v}^{2}-64512\,(\cos(v))^{3}+24192\,(\cos(v))^{2}\\
-22026\,(\cos(v))^{2}{v}^{2}+9656\,\cos(v){v}^{2}+12096\,\cos(v)-2529\,{v}^{2}-4032+64512\,(\cos(v))^{5}
\\
b_{3,num}^0=56448\,(\cos(v))^{4}-73215\,(\cos(v))^{4}{v}^{2}+112896\,(\cos(v))^{3}-23113\,(\cos(v))^{3}{v}^{2}\\
-42336\,(\cos(v))^{2}+73215\,(\cos(v))^{2}{v}^{2}-40011\,\cos(v){v}^{2}-21168\,\cos(v)+10204\,{v}^{2}\\
+7056-112896\,(\cos(v))^{5}
\\
b_{4,num}^0=-225792\,(\cos(v))^{4}+325629\,(\cos(v))^{4}{v}^{2}-451584\,(\cos(v))^{3}-38624\,(\cos(v))^{3}{v}^{2}\\
-96246\,(\cos(v))^{2}{v}^{2}+169344\,(\cos(v))^{2}+28968\,\cos(v){v}^{2}+84672\,\cos(v)+451584\,(\cos(v))^{5}\\
-8047\,{v}^{2}-28224
\\
b_{5,num}^0=-388196\,(\cos(v))^{4}{v}^{2}+282240\,(\cos(v))^{4}+564480\,(\cos(v))^{3}-27081\,(\cos(v))^{3}{v}^{2}\\
+233349\,(\cos(v))^{2}{v}^{2}-211680\,(\cos(v))^{2}-111571\,\cos(v){v}^{2}-105840\,\cos(v)+28899\,{v}^{2}\\
-564480\,(\cos(v))^{5}+35280
\end{array}
\end{equation*}

Since for small values of $v$, the above formulae are subject to
heavy cancelations, we give the Taylor expansions of $b$
coefficients:

\ason
\begin{equation*}
\begin{array}{l}
\ds b_1 = {\frac {399187}{241920}}-{\frac
{52559}{912384}}\,{v}^{2}+{\frac {
100673687}{29059430400}}\,{v}^{4}-{\frac
{1084493}{27897053184}}\,{v}^ {6}+{\frac
{96453547}{213412456857600}}\,{v}^{8}
\\
\ds b_2 = -{\frac {17327}{8640}}+{\frac {52559}{114048}}\,{v}^{2}-{\frac {
100673687}{3632428800}}\,{v}^{4}+{\frac {1084493}{3487131648}}\,{v}^{6
}-{\frac {96453547}{26676557107200}}\,{v}^{8}+\ldots
\\
\ds b_3 = {\frac {597859}{60480}}-{\frac {367913}{228096}}\,{v}^{2}+{\frac {
100673687}{1037836800}}\,{v}^{4}-{\frac {1084493}{996323328}}\,{v}^{6}
+{\frac {96453547}{7621873459200}}\,{v}^{8}-\ldots
\\
\ds b_4= -{\frac {704183}{60480}}+{\frac {367913}{114048}}\,{v}^{2}-{\frac {
100673687}{518918400}}\,{v}^{4}+{\frac {1084493}{498161664}}\,{v}^{6}-
{\frac {96453547}{3810936729600}}\,{v}^{8}+\ldots
\\
\ds b_5 = {\frac {465133}{24192}}-{\frac {1839565}{456192}}\,{v}^{2}+{\frac {
100673687}{415134720}}\,{v}^{4}-{\frac {5422465}{1992646656}}\,{v}^{6}
+{\frac {96453547}{3048749383680}}\,{v}^{8}-\ldots
\end{array}
\end{equation*}
\asoff

The PLTE of the method is given
\begin{equation}
 plte=\left( {\frac {52559}{912384}}\,{w}^{2}y^{(10)}+{\frac {52559}{912384
}}\,y^{(12)} \right) {h}^{12}
\end{equation}
\subsection{Method PF-D1: Phase fitted + 1st Derivative of Phase-lag is zero}
The second method of the family (PF-D1) is constructed by forcing the phase-lag function and its 1st derivative to be zero at the frequency $V=\omega *h$. The coefficients $a$ are left the same, while coefficients $b$ become:

\begin{equation*}
\begin{array}{rlrlrl}
b_0=&0 & b_1=&\frac{1}{192}\frac{b_{1,num}^1}{D_1} & b_2=&\frac{1}{48}\frac{b_{2,num}^1}{D_1} \\
b_3=&\frac{1}{48}\frac{b_{3,num}^1}{D_1} & b_4=&\frac{1}{48}\frac{b_{4,num}^1}{D_1} & b_5=&\frac{1}{96}\frac{b_{5,num}^1}{D_1}
\end{array}
\end{equation*}

where

\begin{equation}
D_1={v}^{3} \left( \left( \cos \left( v \right) \right) ^{5}-3\,
\left( \cos \left( v \right) \right) ^{4}+2\, \left( \cos \left( v
\right)
 \right) ^{3}+2\, \left( \cos \left( v \right) \right) ^{2}-3\,\cos
 \left( v \right) +w \right)
\nonumber
\end{equation}

and

\begin{equation*}
\begin{array}{l}
b_{1,num}^1=-48\,\sin(v)-768\,(\cos(v))^{3}\sin(v)+144\,\cos(v)\sin(v)+768\,\sin(v)(\cos(v))^{5}\\
-384\,(\cos(v))^{4}\sin(v)+288\,(\cos(v))^{2}\sin(v)+432\,\cos(v)v-941\,\cos(v){v}^{3}+281\,{v}^{3}\\
+1344\,(\cos(v))^{5}v+1344\,(\cos(v))^{4}v-1776\,(\cos(v))^{3}v-768\,(\cos(v))^{6}v-576\,(\cos(v))^{2}v\\
+259\,(\cos(v))^{2}{v}^{3}+1481\,(\cos(v))^{3}{v}^{3}
\ec

\bc
b_{2,num}^1=72\,\sin(v)-384\,\sin(v)(\cos(v))^{6}-960\,\sin(v)(\cos(v))^{5}+960\,(\cos(v))^{4}\sin(v)\\
+1008\,(\cos(v))^{3}\sin(v)-504\,(\cos(v))^{2}\sin(v)-192\,\cos(v)\sin(v)-1305\,(\cos(v))^{3}{v}^{3}\\
+732\,(\cos(v))^{2}v+427\,(\cos(v))^{2}{v}^{3}+225\,\cos(v){v}^{3}-612\,\cos(v)v-12\,v-26\,{v}^{3}\\
+192\,(\cos(v))^{7}v+576\,(\cos(v))^{6}v-2304\,(\cos(v))^{5}v-1296\,(\cos(v))^{4}v-1481\,{v}^{3}(\cos(v))^{4}\\
+2724\,(\cos(v))^{3}v\nonumber\\
b_{3,num}^1=-192\,\sin(v)+2304\,\sin(v)(\cos(v))^{6}+1920\,\sin(v)(\cos(v))^{5}\\
-3840\,(\cos(v))^{4}\sin(v)-2208\,(\cos(v))^{3}\sin(v)+1584\,(\cos(v))^{2}\sin(v)+432\,\cos(v)\sin(v)\\
-7464\,(\cos(v))^{3}v+3387\,(\cos(v))^{3}{v}^{3}-1512\,(\cos(v))^{2}v-895\,(\cos(v))^{2}{v}^{3}\\
-1088\,\cos(v){v}^{3}+1512\,\cos(v)v+72\,v-1152\,(\cos(v))^{7}v+384\,(\cos(v))^{6}v+7104\,(\cos(v))^{5}v\\
+1481\,{v}^{3}(\cos(v))^{5}+4443\,{v}^{3}(\cos(v))^{4}+1056\,(\cos
(v))^{4}v+232\,{v}^{3}
\nonumber\\
b_{4,num}^1=-192\,\sin(v)+2304\,\sin(v)(\cos(v))^{6}+1920\,\sin(v)(\cos(v))^{5}\\
-3840\,(\cos(v))^{4}\sin(v)-2208\,(\cos(v))^{3}\sin(v)+1584\,(\cos(v))^{2}\sin(v)+432\,\cos(v)\sin(v)\\
-7464\,(\cos(v))^{3}v+3387\,(\cos(v))^{3}{v}^{3}-1512\,(\cos(v))^{2}v-895\,(\cos(v))^{2}{v}^{3}\\
-1088\,\cos(v){v}^{3}+1512\,\cos(v)v+72\,v-1152\,(\cos(v))^{7}v+384\,(\cos(v))^{6}v+7104\,(\cos(v))^{5}v\\
+1481\,{v}^{3}(\cos(v))^{5}+4443\,{v}^{3}(\cos(v))^{4}+1056\,(\cos
(v))^{4}v+232\,{v}^{3}
\nonumber\\
b_{5,num}^1=-720\,\sin(v)+15360\,\sin(v)(\cos(v))^{6}+3840\,\sin(v)(\cos(v))^{5}\\
-21120\,(\cos(v))^{4}\sin(v)-5760\,(\cos(v))^{3}\sin(v)+7200\,(\cos(v))^{2}\sin(v)+1200\,\cos(v)\sin(v)\\
-29040\,(\cos(v))^{3}v+10059\,(\cos(v))^{3}{v}^{3}-3360\,(\cos(v))^{2}v-4231\,(\cos(v))^{2}{v}^{3}\\
+5040\,\cos(v)v-3411\,\cos(v){v}^{3}+480\,v-7680\,(\cos(v))^{7}v+11520\,(\cos(v))^{6}v\\
+31680\,(\cos(v))^{5}v+12252\,{v}^{3}(\cos(v))^{5}-8640\,(\cos(v))^{4}v+22484\,{v}^{3}(\cos(v))^{4}+647\,{v}^{3}
\end{array}
\end{equation*}

Since for small values of $v$, the above formulae are subject to
heavy cancelations, we give the Taylor expansions of the
coefficients $b$:

\ason
\begin{equation*}
\begin{array}{l}
\ds b_1 = {\frac {399187}{241920}}-{\frac
{52559}{456192}}\,{v}^{2}+{\frac {
975124291}{174356582400}}\,{v}^{4}-{\frac
{2896813}{49816166400}}\,{v} ^{6}+{\frac
{1818828019}{1067062284288000}}\,{v}^{8}
\\
\ds b_2= -{\frac {17327}{8640}}+{\frac {52559}{57024}}\,{v}^{2}-{\frac {
4461254807}{43589145600}}\,{v}^{4}+{\frac {2517959}{340540200}}\,{v}^{
6}-{\frac {98779707713}{266765571072000}}\,{v}^{8}+\ldots
\\
\ds b_3= {\frac {597859}{60480}}-{\frac {367913}{114048}}\,{v}^{2}+{\frac {
3127415341}{6227020800}}\,{v}^{4}-{\frac {3766196569}{87178291200}}\,{
v}^{6}+{\frac {44891085091}{20520428544000}}\,{v}^{8}-\ldots
\\
\ds b_4= -{\frac {704183}{60480}}+{\frac {367913}{57024}}\,{v}^{2}-{\frac {
7330976207}{6227020800}}\,{v}^{4}+{\frac {584032469}{5448643200}}\,{v}
^{6}-{\frac {1452594367391}{266765571072000}}\,{v}^{8}
\\
\ds b_5= {\frac {465133}{24192}}-{\frac {1839565}{228096}}\,{v}^{2}+{\frac {
3844845691}{2490808320}}\,{v}^{4}-{\frac {4974280813}{34871316480}}\,{
v}^{6}+{\frac {773868209533}{106706228428800}}\,{v}^{8}
\end{array}
\end{equation*}
\asoff

The PLTE of the method is given

\begin{equation}
plte=\left( {\frac {52559}{456192}}\,{w}^{2}y^{(10)}+{\frac
{52559}{912384 }}\,{w}^{4}y^{(8)}+{\frac {52559}{912384}}\,y^{(12)}
\right) {h}^{12}
\end{equation}

\subsection{Method PF-D2: Phase fitted + 1st, 2nd Derivatives of Phase-lag are zero}
The third method of the family (PF-D2) is constructed by forcing the phase-lag function and its 1st and 2nd derivatives to be zero at the frequency $V=\omega *h$. The coefficients $a$ are left the same, while for the coefficients $b$ we have:

\begin{equation*}
\begin{array}{rlrlrl}
b_0=&0 \nonumber & b_1=&\frac{1}{32}\frac{b_{1,num}^2}{D_2'}& b_2=&\frac{1}{16}\frac{b_{2,num}^2}{D_2}\\
b_3=&\frac{1}{8}\frac{b_{3,num}^2}{D_2'}& b_4=&\frac{1}{16}\frac{b_{4,num}^2}{D_2}& b_5=&\frac{1}{16}\frac{b_{5,num}^2}{D_2'}
\end{array}
\end{equation*}

where

\begin{equation*}
\begin{array}{l}
D_2= {v}^{4} \left( \left( \cos \left( v \right) \right) ^{5}-3\, \left(
\cos \left( v \right) \right) ^{4}+2\, \left( \cos \left( v \right)
 \right) ^{3}+2\, \left( \cos \left( v \right) \right) ^{2}-3\,\cos
 \left( v \right) +1 \right)
\\
D_2'= {v}^{4} \left( \sin \left( v \right) \right) ^{4} \left( \left( \cos
 \left( v \right) \right) ^{2}-2\,\cos \left( v \right) +1 \right)
\end{array}
\end{equation*}

and

\begin{equation*}
\begin{array}{l}
b_{1,num}^2 = -48\,\sin(v)-768\,(\cos(v))^{3}\sin(v)+144\,\cos(v)\sin(v)+768\,\sin(v)(\cos(v))^{5}\\
-384\,(\cos(v))^{4}\sin(v)+288\,(\cos(v))^{2}\sin(v)+432\,\cos(v)v-941\,\cos(v){v}^{3}+281\,{v}^{3}\\
+1344\,(\cos(v))^{5}v+1344\,(\cos(v))^{4}v-1776\,(\cos(v))^{3}v-768\,(\cos(v))^{6}v-576\,(\cos(v))^{2}v\\
+259\,(\cos(v))^{2}{v}^{3}+1481\,(\cos(v))^{3}{v}^{3}
\\
b_{2,num}^2 = -24-800\,v(\cos(v))^{3}\sin(v)-512\,\sin(v)v(\cos(v))^{6}+1088\,v(\cos(v))^{4}\sin(v)\\
-25\,{v}^{4}+628\,(\cos(v))^{3}{v}^{2}-644\,(\cos(v))^{2}{v}^{2}+216\,v\cos(v)\sin(v)+24\,v\sin(v)\\
+168\,(\cos(v))^{2}-528\,v(\cos(v))^{2}\sin(v)-640\,(\cos(v))^{5}{v}^{2}+68\,{v}^{2}+768\,(\cos(v))^{5}\\
+192\,(\cos(v))^{6}-180\,\cos(v){v}^{2}+385\,(\cos(v))^{2}{v}^{4}+72\,\cos(v)-456\,(\cos(v))^{3}-\\
336\,(\cos(v))^{4}+512\,\sin(v)v(\cos(v))^{5}+435\,(\cos(v))^{3}{v}^{4}+768\,(\cos(v))^{4}{v}^{2}\\
-384\,(\cos(v))^{7}+192\,{v}^{2}(\cos(v))^{7}-192\,{v}^{2}(\cos(v))^{6}-75\,\cos(v){v}^{4}
\\
b_{3,num}^2 = 24+1612\,v(\cos(v))^{3}\sin(v)+736\,\sin(v)v(\cos(v))^{7}-960\,\sin(v)v(\cos(v))^{6}\\
+676\,v(\cos(v))^{4}\sin(v)+55\,{v}^{4}+548\,(\cos(v))^{3}{v}^{2}+804\,(\cos(v))^{2}{v}^{2}\\
-240\,v\cos(v)\sin(v)-32\,v\sin(v)+192\,\sin(v)(\cos(v))^{8}v-1050\,(\cos(v))^{4}{v}^{4}\\
+672\,(\cos(v))^{8}-300\,(\cos(v))^{2}+136\,v(\cos(v))^{2}\sin(v)-1028\,(\cos(v))^{5}{v}^{2}\\
-60\,{v}^{2}+132\,(\cos(v))^{5}-1560\,(\cos(v))^{6}-435\,{v}^{4}(\cos(v))^{5}+96\,\cos(v){v}^{2}\\
-256\,(\cos(v))^{8}{v}^{2}-265\,(\cos(v))^{2}{v}^{4}-24\,\cos(v)+84\,(\cos(v))^{3}+1164\,(\cos(v))^{4}\\
-2120\,\sin(v)v(\cos(v))^{5}-865\,(\cos(v))^{3}{v}^{4}-1264\,(\cos(v))^{4}{v}^{2}-64\,(\cos(v))^{9}{v}^{2}\\
-384\,(\cos(v))^{7}+192\,(\cos(v))^{9}+448\,{v}^{2}(\cos(v))^{7}+776\,{v}^{2}(\cos(v))^{6}+40\,\cos(v){v}^{4}
\ec

\bc
b_{4,num}^2 = -72-3904\,v(\cos(v))^{3}\sin(v)-1536\,\sin(v)v(\cos(v))^{7}+256\,\sin(v)v(\cos(v))^{6}\\
+512\,v(\cos(v))^{4}\sin(v)-75\,{v}^{4}-244\,(\cos(v))^{3}{v}^{2}-1980\,(\cos(v))^{2}{v}^{2}\\
+648\,v\cos(v)\sin(v)+104\,v\sin(v)+1740\,(\cos(v))^{4}{v}^{4}-1536\,(\cos(v))^{8}+792\,(\cos(v))^{2}\\
-944\,v(\cos(v))^{2}\sin(v)+832\,(\cos(v))^{5}{v}^{2}+156\,{v}^{2}+960\,(\cos(v))^{5}+3648\,(\cos(v))^{6}\\
+580\,{v}^{4}(\cos(v))^{5}-396\,\cos(v){v}^{2}+512\,(\cos(v))^{8}{v}^{2}+855\,(\cos(v))^{2}{v}^{4}+120\,\cos(v)\\
-696\,(\cos(v))^{3}-2832\,(\cos(v))^{4}+4864\,\sin(v)v(\cos(v))^{5}+2165\,(\cos(v))^{3}{v}^{4}\\
+3168\,(\cos(v))^{4}{v}^{2}-384\,(\cos(v))^{7}-192\,{v}^{2}(\cos(v))^{7}-1856\,{v}^{2}(\cos(v))^{6}-225\,\cos(v){v}^{4}
\\
b_{5,num}^2 = 84+5416\,v(\cos(v))^{3}\sin(v)+640\,\sin(v)v(\cos(v))^{7}-6720\,\sin(v)v(\cos(v))^{6}\\
+4080\,v(\cos(v))^{4}\sin(v)+165\,{v}^{4}+2540\,(\cos(v))^{3}{v}^{2}+2764\,(\cos(v))^{2}{v}^{2}\\
-900\,v\cos(v)\sin(v)-120\,v\sin(v)+2304\,\sin(v)(\cos(v))^{8}v-4720\,(\cos(v))^{4}{v}^{4}\\
+1920\,(\cos(v))^{8}-1020\,(\cos(v))^{2}+420\,v(\cos(v))^{2}\sin(v)-920\,{v}^{4}(\cos(v))^{6}\\
-4528\,(\cos(v))^{5}{v}^{2}-172\,{v}^{2}+2736\,(\cos(v))^{5}-4800\,(\cos(v))^{6}-3380\,{v}^{4}(\cos(v))^{5}\\
+260\,\cos(v){v}^{2}-825\,(\cos(v))^{2}{v}^{4}-60\,\cos(v)-180\,(\cos(v))^{3}+3816\,(\cos(v))^{4}\\
-5120\,\sin(v)v(\cos(v))^{5}-3115\,(\cos(v))^{3}{v}^{4}-2912\,(\cos(v))^{4}{v}^{2}-768\,(\cos(v))^{9}{v}^{2}\\
-4800\,(\cos(v))^{7}+2304\,(\cos(v))^{9}+2496\,{v}^{2}(\cos(v))^{7}+320\,{v}^{2}(\cos(v))^{6}+195\,\cos(v){v}^{4}
\end{array}
\end{equation*}

Since for small values of $v$, the above formulae are subject to
heavy cancelations, we give the Taylor expansions of the
coefficients $b$:

\ason
\begin{equation*}
\begin{array}{l}
\ds b_1 = {\frac {399187}{241920}}-{\frac
{52559}{304128}}\,{v}^{2}+{\frac {
371082169}{58118860800}}\,{v}^{4}-{\frac
{83360891}{523069747200}}\,{v }^{6}-{\frac
{1467578899}{355687428096000}}\,{v}^{8}
\\
\ds b_2= -{\frac {17327}{8640}}+{\frac {52559}{38016}}\,{v}^{2}-{\frac {
3253170563}{14529715200}}\,{v}^{4}+{\frac {5070942803}{261534873600}}
\,{v}^{6}-{\frac {86978398867}{88921857024000}}\,{v}^{8}
\\
\ds b_3={\frac {597859}{60480}}-{\frac {367913}{76032}}\,{v}^{2}+{\frac {
2523373219}{2075673600}}\,{v}^{4}-{\frac {22329042629}{130767436800}}
\,{v}^{6}+{\frac {1453392734357}{88921857024000}}\,{v}^{8}
\ec

\bc
\ds b_4= -{\frac {704183}{60480}}+{\frac {367913}{38016}}\,{v}^{2}-{\frac {
6122891963}{2075673600}}\,{v}^{4}+{\frac {133660742933}{261534873600}}
\,{v}^{6}-{\frac {5024895032029}{88921857024000}}\,{v}^{8}
\\
\ds b_5= {\frac {465133}{24192}}-{\frac {1839565}{152064}}\,{v}^{2}+{\frac {
3240803569}{830269440}}\,{v}^{4}-{\frac {37612768013}{52306974720}}\,{
v}^{6}+{\frac {2927078073011}{35568742809600}}\,{v}^{8}
\end{array}
\end{equation*}
\asoff

The PLTE of the method is given

\begin{equation}
plte=\left( {\frac {52559}{912384}}\,{w}^{6}y^{(6)}+{\frac
{52559}{304128} }\,{w}^{2}y^{(10)}+{\frac
{52559}{912384}}\,y^{(12)}+{\frac {52559}{ 304128}}\,{w}^{4}y^{(8)}
\right) {h}^{12}
\end{equation}

\subsection{Method PF-D3: Phase fitted + 1st, 2nd, 3rd Derivatives of Phase-lag are zero}
The fourth method of the family (PF-D3) is constructed by forcing the phase-lag function and its 1st, 2nd and 3rd derivatives to be zero at the frequency $V=\omega *h$. The coefficients $a$ are left the same, while for the coefficients $b$ we have:

\begin{equation*}
\begin{array}{rlrlrl}
b_0=&0 \nonumber & b_1=&\frac{1}{48}\frac{b_{1,num}^3}{D_3} & b_2=&\frac{1}{12}\frac{b_{2,num}^3}{D_3} \\
b_3=&\frac{1}{12}\frac{b_{3,num}^3}{D_3} & b_4=&\frac{1}{12}\frac{b_{4,num}^3}{D_3} & b_5=&\frac{1}{24}\frac{b_{5,num}^3}{D_3}
\end{array}
\end{equation*}

where

\begin{equation*}
\begin{array}{l}
D_3={v}^{5}((\cos(v))^{7}-(\cos(v))^{6}-3\,(\cos(v))^{5}+3\,(\cos(v))^{4}+3\,(\cos(v))^{3}-3\,(\cos(v))^{2}\\
-\cos(v)+1)
\end{array}
\end{equation*}

and

\begin{equation*}
\begin{array}{l}
b_{1,num}^3=-2096\,{v}^{2}(\cos(v))^{5}\sin(v)+1508\,{v}^{2}(\cos(v))^{4}\sin(v)+1648\,{v}^{2}(\cos(v))^{3}\sin(v)\\
-964\,{v}^{2}\sin(v)(\cos(v))^{2}+36\,v-366\,{v}^{2}\cos(v)\sin(v)+135\,{v}^{5}\cos(v)-48\,{v}^{3}-\\
384\,(\cos(v))^{8}{v}^{3}-336\,(\cos(v))^{4}\sin(v)-656\,{v}^{2}(\cos(v))^{6}\sin(v)-336\,\cos(v){v}^{3}\\
+1440\,(\cos(v))^{5}v-738\,(\cos(v))^{2}v-384\,\sin(v)(\cos(v))^{7}+168\,(\cos(v))^{2}\sin(v)\\
+234\,\cos(v)v+528\,(\cos(v))^{2}{v}^{3}-1098\,(\cos(v))^{3}v-2448\,(\cos(v))^{6}v\\
+768\,\sin(v)(\cos(v))^{5}+192\,\sin(v)(\cos(v))^{6}-456\,(\cos(v))^{3}\sin(v)+1008\,(\cos(v))^{3}{v}^{3}\\
+72\,\cos(v)\sin(v)+864\,v(\cos(v))^{8}+2286\,(\cos(v))^{4}v+45\,{v}^{5}(\cos(v))^{3}-24\,\sin(v)\\
+1200\,(\cos(v))^{6}{v}^{3}-1008\,(\cos(v))^{5}{v}^{3}+45\,{v}^{5}+135\,{v}^{5}(\cos(v))^{2}\\
-576\,(\cos(v))^{7}v-1296\,(\cos(v))^{4}{v}^{3}+832\,\sin(v){v}^{2}(\cos(v))^{7}+94\,{v}^{2}\sin(v)\\
+336\,(\cos(v))^{7}{v}^{3}
\\
b_{2,num}^3=y-692\,{v}^{2}(\cos(v))^{5}\sin(v)-2407\,{v}^{2}(\cos(v))^{4}\sin(v)+235\,{v}^{2}(\cos(v))^{3}\sin(v)\\
+776\,{v}^{2}\sin(v)(\cos(v))^{2}-27\,v+123\,{v}^{2}\cos(v)\sin(v)-90\,{v}^{5}\cos(v)+48\,{v}^{3}\\
-96\,(\cos(v))^{8}{v}^{3}+852\,(\cos(v))^{4}\sin(v)+2560\,{v}^{2}(\cos(v))^{6}\sin(v)+144\,\cos(v){v}^{3}\\
-3690\,(\cos(v))^{5}v+90\,(\cos(v))^{2}v-96\,\sin(v)(\cos(v))^{7}-192\,(\cos(v))^{2}\sin(v)-144\,\cos(v)v\\
-48\,(\cos(v))^{2}{v}^{3}+1458\,(\cos(v))^{3}v-738\,(\cos(v))^{6}v+120\,\sin(v)(\cos(v))^{5}\\
-1248\,\sin(v)(\cos(v))^{6}-24\,(\cos(v))^{3}\sin(v)-912\,\sin(v){v}^{2}(\cos(v))^{8}-720\,(\cos(v))^{3}{v}^{3}\\
+360\,v(\cos(v))^{8}+315\,(\cos(v))^{4}v+576\,\sin(v)(\cos(v))^{8}-270\,{v}^{5}(\cos(v))^{3}-1152\,(\cos(v))^{9}v\\
+12\,\sin(v)+240\,(\cos(v))^{6}{v}^{3}+1296\,(\cos(v))^{5}{v}^{3}-270\,{v}^{5}(\cos(v))^{2}-90\,(\cos(v))^{4}{v}^{5}\\
+3528\,(\cos(v))^{7}v-144\,(\cos(v))^{4}{v}^{3}+352\,\sin(v){v}^{2}(\cos(v))^{7}+288\,(\cos(v))^{9}{v}^{3}\\
-35\,{v}^{2}\sin(v)-1008\,(\cos(v))^{7}{v}^{3}
\ec

\bc
b_{3,num}^3=3664\,{v}^{2}(\cos(v))^{5}\sin(v)+3218\,{v}^{2}(\cos(v))^{4}\sin(v)-812\,{v}^{2}(\cos(v))^{3}\sin(v)\\
-1498\,{v}^{2}\sin(v)(\cos(v))^{2}+54\,v-258\,{v}^{2}\cos(v)\sin(v)+135\,{v}^{5}\cos(v)-84\,{v}^{3}\\
+1344\,(\cos(v))^{8}{v}^{3}-1200\,(\cos(v))^{4}\sin(v)-2420\,{v}^{2}(\cos(v))^{6}\sin(v)\\
+1344\,\sin(v){v}^{2}(\cos(v))^{9}-300\,\cos(v){v}^{3}+4842\,(\cos(v))^{5}v-180\,(\cos(v))^{2}v\\
+2496\,\sin(v)(\cos(v))^{7}-384\,(\cos(v))^{10}{v}^{3}+312\,(\cos(v))^{2}\sin(v)+324\,\cos(v)v\\
+60\,(\cos(v))^{2}{v}^{3}-2502\,(\cos(v))^{3}v+6498\,(\cos(v))^{6}v-1608\,\sin(v)(\cos(v))^{5}\\
+270\,(\cos(v))^{5}{v}^{5}+1488\,\sin(v)(\cos(v))^{6}+264\,(\cos(v))^{3}\sin(v)+624\,\sin(v){v}^{2}(\cos(v))^{8}\\
+1092\,(\cos(v))^{3}{v}^{3}-6336\,v(\cos(v))^{8}-2052\,(\cos(v))^{4}v-576\,\sin(v)(\cos(v))^{8}\\
+855\,{v}^{5}(\cos(v))^{3}+864\,(\cos(v))^{9}v-24\,\sin(v)-1644\,(\cos(v))^{6}{v}^{3}-1476\,(\cos(v))^{5}{v}^{3}\\
+45\,{v}^{5}+405\,{v}^{5}(\cos(v))^{2}-1152\,\sin(v)(\cos(v))^{9}+810\,(\cos(v))^{4}{v}^{5}-3528\,(\cos(v))^{7}v\\
+708\,(\cos(v))^{4}{v}^{3}-3920\,\sin(v){v}^{2}(\cos(v))^{7}-192\,(\cos(v))^{9}{v}^{3}+2016\,v(\cos(v))^{10}\\
+58\,{v}^{2}\sin(v)+876\,(\cos(v))^{7}{v}^{3}
\\
b_{4,num}^3=-5752\,{v}^{2}(\cos(v))^{5}\sin(v)-3789\,{v}^{2}(\cos(v))^{4}\sin(v)+1333\,{v}^{2}(\cos(v))^{3}\sin(v)\\
+2196\,{v}^{2}\sin(v)(\cos(v))^{2}-81\,v+417\,{v}^{2}\cos(v)\sin(v)-270\,{v}^{5}\cos(v)+108\,{v}^{3}\\
-2016\,(\cos(v))^{8}{v}^{3}+1404\,(\cos(v))^{4}\sin(v)+1036\,{v}^{2}(\cos(v))^{6}\sin(v)\\
-2048\,\sin(v){v}^{2}(\cos(v))^{9}+468\,\cos(v){v}^{3}-5310\,(\cos(v))^{5}v+234\,(\cos(v))^{2}v\\
-4224\,\sin(v)(\cos(v))^{7}+576\,(\cos(v))^{10}{v}^{3}-432\,(\cos(v))^{2}\sin(v)-504\,\cos(v)v\\
-36\,(\cos(v))^{2}{v}^{3}+3582\,(\cos(v))^{3}v-10818\,(\cos(v))^{6}v+2760\,\sin(v)(\cos(v))^{5}\\
-1080\,(\cos(v))^{5}{v}^{5}-816\,\sin(v)(\cos(v))^{6}-456\,(\cos(v))^{3}\sin(v)+1360\,\sin(v){v}^{2}(\cos(v))^{8}\\
-1500\,(\cos(v))^{3}{v}^{3}-1152\,(\cos(v))^{11}v+10224\,v(\cos(v))^{8}+3609\,(\cos(v))^{4}v\\
-960\,\sin(v)(\cos(v))^{8}-1170\,{v}^{5}(\cos(v))^{3}+2592\,(\cos(v))^{9}v+36\,\sin(v)\\
+768\,\sin(v)(\cos(v))^{10}+2484\,(\cos(v))^{6}{v}^{3}+1500\,(\cos(v))^{5}{v}^{3}+192\,(\cos(v))^{11}{v}^{3}\\
-810\,{v}^{5}(\cos(v))^{2}-704\,\sin(v){v}^{2}(\cos(v))^{10}+1920\,\sin(v)(\cos(v))^{9}-1350\,(\cos(v))^{4}{v}^{5}\\
+792\,(\cos(v))^{7}v-1116\,(\cos(v))^{4}{v}^{3}+6032\,\sin(v){v}^{2}(\cos(v))^{7}-360\,{v}^{5}(\cos(v))^{6}\\
-480\,(\cos(v))^{9}{v}^{3}-3168\,v(\cos(v))^{10}-81\,{v}^{2}\sin(v)-180\,(\cos(v))^{7}{v}^{3}
\ec

\bc
b_{5,num}^3=13216\,{v}^{2}(\cos(v))^{5}\sin(v)+10404\,{v}^{2}(\cos(v))^{4}\sin(v)-4672\,{v}^{2}(\cos(v))^{3}\sin(v)\\
-4932\,{v}^{2}\sin(v)(\cos(v))^{2}+180\,v-762\,{v}^{2}\cos(v)\sin(v)+405\,{v}^{5}\cos(v)-240\,{v}^{3}\\
+3456\,(\cos(v))^{8}{v}^{3}-3888\,(\cos(v))^{4}\sin(v)-4048\,{v}^{2}(\cos(v))^{6}\sin(v)\\
+2816\,\sin(v){v}^{2}(\cos(v))^{9}-912\,\cos(v){v}^{3}+15192\,(\cos(v))^{5}v+162\,(\cos(v))^{2}v\\
+7680\,\sin(v)(\cos(v))^{7}-768\,(\cos(v))^{10}{v}^{3}+1080\,(\cos(v))^{2}\sin(v)+1062\,\cos(v)v\\
-432\,(\cos(v))^{2}{v}^{3}-9054\,(\cos(v))^{3}v+22680\,(\cos(v))^{6}v-5856\,\sin(v)(\cos(v))^{5}\\
+2160\,(\cos(v))^{5}{v}^{5}+2112\,\sin(v)(\cos(v))^{6}+1320\,(\cos(v))^{3}\sin(v)-4288\,\sin(v){v}^{2}(\cos(v))^{8}\\
+3504\,(\cos(v))^{3}{v}^{3}+4608\,(\cos(v))^{11}v-72\,\cos(v)\sin(v)-17856\,v(\cos(v))^{8}\\
-9774\,(\cos(v))^{4}v+3840\,\sin(v)(\cos(v))^{8}+3375\,{v}^{5}(\cos(v))^{3}-9216\,(\cos(v))^{9}v\\
-72\,\sin(v)-3072\,\sin(v)(\cos(v))^{10}-5520\,(\cos(v))^{6}{v}^{3}-4272\,(\cos(v))^{5}{v}^{3}+135\,{v}^{5}\\
-768\,(\cos(v))^{11}{v}^{3}+1485\,{v}^{5}(\cos(v))^{2}+2816\,\sin(v){v}^{2}(\cos(v))^{10}\\
-3072\,\sin(v)(\cos(v))^{9}+3600\,(\cos(v))^{4}{v}^{5}-2592\,(\cos(v))^{7}v+3504\,(\cos(v))^{4}{v}^{3}\\
-10688\,\sin(v){v}^{2}(\cos(v))^{7}+360\,{v}^{5}(\cos(v))^{7}+1080\,{v}^{5}(\cos(v))^{6}\\
+1536\,(\cos(v))^{9}{v}^{3}+4608\,v(\cos(v))^{10}+138\,{v}^{2}\sin(v)+912\,(\cos(v))^{7}{v}^{3}
\end{array}
\end{equation*}

Since for small values of $v$, the above formulae are subject to
heavy cancelations, we give the Taylor expansions of the
coefficients $b$:

\ason
\begin{equation*}
\begin{array}{l}
\ds b_1 = {\frac {399187}{241920}}-{\frac
{52559}{228096}}\,{v}^{2}+{\frac {
11315653}{1937295360}}\,{v}^{4}-{\frac
{5807033}{13076743680}}\,{v}^{6 }-{\frac
{614853845}{17072996548608}}\,{v}^{8}
 \\
\ds b_2= -{\frac {17327}{8640}}+{\frac {52559}{28512}}\,{v}^{2}-{\frac {5758537
}{14676480}}\,{v}^{4}+{\frac {225159101}{6538371840}}\,{v}^{6}-{\frac
{8200289261}{4268249137152}}\,{v}^{8}
 \\
\ds b_3= {\frac {597859}{60480}}-{\frac {367913}{57024}}\,{v}^{2}+{\frac {
154801723}{69189120}}\,{v}^{4}-{\frac {2799488011}{6538371840}}\,{v}^{
6}+{\frac {1063054198007}{21341245685760}}\,{v}^{8}
 \\
\ds b_4= -{\frac {704183}{60480}}+{\frac {367913}{28512}}\,{v}^{2}-{\frac {
381346481}{69189120}}\,{v}^{4}+{\frac {9217976399}{6538371840}}\,{v}^{
6}-{\frac {5100295346143}{21341245685760}}\,{v}^{8}
 \\
\ds b_5= {\frac {465133}{24192}}-{\frac {1839565}{114048}}\,{v}^{2}+{\frac {
67543471}{9225216}}\,{v}^{4}-{\frac {241481599}{118879488}}\,{v}^{6}+{
\frac {16316044646989}{42682491371520}}\,{v}^{8}
\end{array}
\end{equation*}
\asoff

The PLTE of the method is given

\bc
\ds plte=\left( {\frac {52559}{228096}}\,{w}^{2}y^{(10)}+{\frac
{52559}{912384 }}\,{w}^{8}y^{(4)}+{\frac
{52559}{152064}}\,{w}^{4}y^{(8)}\right.\\
\ds \left.+{\frac {
52559}{228096}}\,{w}^{6}y^{(6)}+{\frac {52559}{912384}}\,y^{(12)}
 \right) {h}^{12}
\ec

\subsection{Method PF-D4: Phase fitted + 1st, 2nd, 3rd, 4th Derivatives of Phase-lag are zero}
The fifth method of the family (PF-D4) is constructed by forcing the phase-lag function and its 1st, 2nd, 3rd and 4th derivatives to be zero at the frequency $V=\omega *h$. The coefficients $a$ are left the same, while for the coefficients $b$ we have:

\begin{equation*}
\begin{array}{rlrlrl}
b_0=&0 \nonumber & b_1=&\frac{1}{96}\frac{b_{1,num}^4}{D_4} & b_2=&\frac{1}{12}\frac{b_{2,num}^4}{D_4} \\
b_3=&\frac{1}{24}\frac{b_{3,num}^4}{D_4} & b_4=&\frac{1}{12}\frac{b_{4,num}^4}{D_4} & b_5=&\frac{1}{48}\frac{b_{5,num}^4}{D_4}
\end{array}
\end{equation*}

where

\begin{equation*}
\begin{array}{l}
D_4={v}^{6} \left( \cos \left( v \right) +1 \right) \left( \sin \left( v
 \right) \right) ^{5}
\end{array}
\end{equation*}

and

\begin{equation*}
\begin{array}{l}
b_{1,num}^4=3069\,\sin(v)(\cos(v))^{2}{v}^{2}-6408\,\sin(v)(\cos(v))^{4}{v}^{2}-4218\,\sin(v)(\cos(v))^{3}{v}^{2}\\
-864\,\sin(v)\cos(v){v}^{4}+1728\,\sin(v)(\cos(v))^{3}{v}^{4}+2016\,\sin(v)(\cos(v))^{4}{v}^{4}\\
-864\,\sin(v)(\cos(v))^{5}{v}^{4}+1380\,\sin(v)\cos(v){v}^{2}-1152\,\sin(v)(\cos(v))^{2}{v}^{4}\\
+2568\,\sin(v)(\cos(v))^{5}{v}^{2}+3408\,\sin(v)(\cos(v))^{6}{v}^{2}-960\,\sin(v)(\cos(v))^{6}{v}^{4}\\
-144\,v+1440\,(\cos(v))^{4}\sin(v)+1029\,\cos(v){v}^{3}+60\,\sin(v)+372\,{v}^{3}+4836\,(\cos(v))^{4}{v}^{3}\\
-5952\,(\cos(v))^{5}v+1728\,(\cos(v))^{2}v-540\,(\cos(v))^{2}\sin(v)-648\,\cos(v)v-3234\,(\cos(v))^{2}{v}^{3}\\
+3912\,(\cos(v))^{3}v+1728\,(\cos(v))^{6}v-480\,\sin(v)(\cos(v))^{5}+2688\,(\cos(v))^{7}v\\
-2064\,(\cos(v))^{6}{v}^{3}+600\,(\cos(v))^{3}\sin(v)-204\,\sin(v){v}^{2}-4786\,(\cos(v))^{3}{v}^{3}\\
+96\,\sin(v){v}^{4}-120\,\cos(v)\sin(v)-3312\,(\cos(v))^{4}v+6176\,(\cos(v))^{5}{v}^{3}\\
-960\,\sin(v)(\cos(v))^{6}-2464\,(\cos(v))^{7}{v}^{3}
\\
b_{2,num}^4=-1047\,\sin(v)(\cos(v))^{2}{v}^{2}+3552\,\sin(v)(\cos(v))^{4}{v}^{2}-2295\,\sin(v)(\cos(v))^{3}{v}^{2}\\
+48\,\sin(v)\cos(v){v}^{4}+384\,\sin(v)(\cos(v))^{3}{v}^{4}-912\,\sin(v)(\cos(v))^{4}{v}^{4}\\
-912\,\sin(v)(\cos(v))^{5}{v}^{4}+78\,\sin(v)\cos(v){v}^{2}+384\,\sin(v)(\cos(v))^{2}{v}^{4}\\
+5184\,\sin(v)(\cos(v))^{5}{v}^{2}-2208\,\sin(v)(\cos(v))^{6}{v}^{2}-2832\,\sin(v){v}^{2}(\cos(v))^{7}\\
+480\,\sin(v)(\cos(v))^{6}{v}^{4}+12\,v-600\,(\cos(v))^{4}\sin(v)-131\,\cos(v){v}^{3}-35\,{v}^{3}\\
+2929\,(\cos(v))^{4}{v}^{3}+3096\,(\cos(v))^{5}v+528\,(\cos(v))^{2}v+120\,(\cos(v))^{2}\sin(v)+120\,\cos(v)v\\
-433\,(\cos(v))^{2}{v}^{3}-1584\,(\cos(v))^{3}v+5472\,(\cos(v))^{6}v-1440\,\sin(v)(\cos(v))^{5}\\
-1632\,(\cos(v))^{7}v-4128\,(\cos(v))^{6}{v}^{3}+540\,(\cos(v))^{3}\sin(v)-27\,\sin(v){v}^{2}\\
+2181\,(\cos(v))^{3}{v}^{3}+48\,\sin(v){v}^{4}-60\,\cos(v)\sin(v)-3516\,(\cos(v))^{4}v-3480\,(\cos(v))^{5}{v}^{3}\\
+480\,\sin(v)(\cos(v))^{6}+1520\,(\cos(v))^{7}{v}^{3}+1712\,(\cos(v))^{8}{v}^{3}\\
-2496\,(\cos(v))^{8}v+960\,\sin(v)(\cos(v))^{7}+480\,\sin(v)(\cos(v))^{7}{v}^{4}
\\
b_{3,num}^4=2997\,\sin(v)(\cos(v))^{2}{v}^{2}+3906\,\sin(v)(\cos(v))^{4}{v}^{2}+498\,\sin(v)(\cos(v))^{3}{v}^{2}\\
-936\,\sin(v)\cos(v){v}^{4}-48\,\sin(v)(\cos(v))^{3}{v}^{4}+24\,\sin(v)(\cos(v))^{4}{v}^{4}\\
+2904\,\sin(v)(\cos(v))^{5}{v}^{4}+1632\,\sin(v)\cos(v){v}^{2}-1008\,\sin(v)(\cos(v))^{2}{v}^{4}\\
-15396\,\sin(v)(\cos(v))^{5}{v}^{2}-21792\,\sin(v)(\cos(v))^{6}{v}^{2}-7488\,(\cos(v))^{9}{v}^{3}
\ec

\bc
+11376\,\sin(v){v}^{2}(\cos(v))^{7}+2880\,\sin(v)(\cos(v))^{6}{v}^{4}-144\,v-1800\,(\cos(v))^{4}\sin(v)\\
+1101\,\cos(v){v}^{3}+60\,\sin(v)+300\,{v}^{3}-4740\,(\cos(v))^{4}{v}^{3}+12768\,(\cos(v))^{5}v\\
+1152\,(\cos(v))^{2}v-180\,(\cos(v))^{2}\sin(v)-792\,\cos(v)v-2910\,(\cos(v))^{2}{v}^{3}\\
+1464\,(\cos(v))^{3}v-15552\,(\cos(v))^{6}v+3120\,\sin(v)(\cos(v))^{5}-27264\,(\cos(v))^{7}v\\
+13632\,(\cos(v))^{6}{v}^{3}-120\,(\cos(v))^{3}\sin(v)-168\,\sin(v){v}^{2}-3244\,(\cos(v))^{3}{v}^{3}\\
+24\,\sin(v){v}^{4}-120\,\cos(v)\sin(v)+5328\,(\cos(v))^{4}v-5980\,(\cos(v))^{5}{v}^{3}\\
+7680\,\sin(v)(\cos(v))^{6}+15296\,(\cos(v))^{7}{v}^{3}-6912\,(\cos(v))^{8}{v}^{3}+9216\,(\cos(v))^{8}v\\
-2880\,\sin(v)(\cos(v))^{7}+14112\,\sin(v)(\cos(v))^{8}{v}^{2}+13824\,(\cos(v))^{9}v\\
-1920\,\sin(v)(\cos(v))^{7}{v}^{4}-1920\,\sin(v)(\cos(v))^{8}{v}^{4}-5760\,\sin(v)(\cos(v))^{8}
\\
b_{4,num}^4=-3321\,\sin(v)(\cos(v))^{2}{v}^{2}+8268\,\sin(v)(\cos(v))^{4}{v}^{2}-6861\,\sin(v)(\cos(v))^{3}{v}^{2}\\
+108\,\sin(v)\cos(v){v}^{4}+1224\,\sin(v)(\cos(v))^{3}{v}^{4}-1812\,\sin(v)(\cos(v))^{4}{v}^{4}\\
-1812\,\sin(v)(\cos(v))^{5}{v}^{4}+162\,\sin(v)\cos(v){v}^{2}+1224\,\sin(v)(\cos(v))^{2}{v}^{4}\\
+10188\,\sin(v)(\cos(v))^{5}{v}^{2}+3552\,\sin(v)(\cos(v))^{6}{v}^{2}+960\,\sin(v)(\cos(v))^{9}{v}^{4}\\
+3648\,(\cos(v))^{9}{v}^{3}+5328\,\sin(v){v}^{2}(\cos(v))^{7}-7872\,\sin(v){v}^{2}(\cos(v))^{9}\\
-480\,\sin(v)(\cos(v))^{6}{v}^{4}+36\,v-1320\,(\cos(v))^{4}\sin(v)-297\,\cos(v){v}^{3}-81\,{v}^{3}\\
+7815\,(\cos(v))^{4}{v}^{3}+4104\,(\cos(v))^{5}v+1728\,(\cos(v))^{2}v+360\,(\cos(v))^{2}\sin(v)\\
+360\,\cos(v)v-1287\,(\cos(v))^{2}{v}^{3}-4464\,(\cos(v))^{3}v+5520\,(\cos(v))^{6}v\\
-2160\,\sin(v)(\cos(v))^{5}+5760\,(\cos(v))^{7}v-5828\,(\cos(v))^{6}{v}^{3}+1380\,(\cos(v))^{3}\sin(v)\\
-81\,\sin(v){v}^{2}+6315\,(\cos(v))^{3}{v}^{3}+108\,\sin(v){v}^{4}-180\,\cos(v)\sin(v)-9396\,(\cos(v))^{4}v\\
-5244\,(\cos(v))^{5}{v}^{3}-960\,\sin(v)(\cos(v))^{6}-3792\,(\cos(v))^{7}{v}^{3}+3840\,\sin(v)(\cos(v))^{9}\\
+3904\,(\cos(v))^{10}{v}^{3}-4208\,(\cos(v))^{8}{v}^{3}+10560\,(\cos(v))^{8}v-2880\,\sin(v)(\cos(v))^{7}\\
-6528\,\sin(v)(\cos(v))^{8}{v}^{2}-5760\,(\cos(v))^{9}v-480\,\sin(v)(\cos(v))^{7}{v}^{4}\\
+960\,\sin(v)(\cos(v))^{8}{v}^{4}+1920\,\sin(v)(\cos(v))^{8}-8448\,(\cos(v))^{10}v
\ec

\bc
b_{5,num}^4=9207\,\sin(v)(\cos(v))^{2}{v}^{2}+21552\,\sin(v)(\cos(v))^{4}{v}^{2}+6498\,\sin(v)(\cos(v))^{3}{v}^{2}\\
-2784\,\sin(v)\cos(v){v}^{4}-2112\,\sin(v)(\cos(v))^{3}{v}^{4}-1824\,\sin(v)(\cos(v))^{4}{v}^{4}\\
+11040\,\sin(v)(\cos(v))^{5}{v}^{4}+5148\,\sin(v)\cos(v){v}^{2}-3072\,\sin(v)(\cos(v))^{2}{v}^{4}\\
-60312\,\sin(v)(\cos(v))^{5}{v}^{2}-81288\,\sin(v)(\cos(v))^{6}{v}^{2}-1536\,\sin(v)(\cos(v))^{9}{v}^{4}\\
-14592\,(\cos(v))^{9}{v}^{3}+27888\,\sin(v){v}^{2}(\cos(v))^{7}+11328\,\sin(v){v}^{2}(\cos(v))^{9}\\
+10944\,\sin(v)(\cos(v))^{6}{v}^{4}-432\,v-8160\,(\cos(v))^{4}\sin(v)+3183\,\cos(v){v}^{3}\\
+180\,\sin(v)+828\,{v}^{3}-23348\,(\cos(v))^{4}{v}^{3}+55680\,(\cos(v))^{5}v+2880\,(\cos(v))^{2}v\\
-180\,(\cos(v))^{2}\sin(v)-2520\,\cos(v)v-8598\,(\cos(v))^{2}{v}^{3}+1560\,(\cos(v))^{3}v\\
-54720\,(\cos(v))^{6}v+12000\,\sin(v)(\cos(v))^{5}-93120\,(\cos(v))^{7}v+47712\,(\cos(v))^{6}{v}^{3}\\
-1080\,(\cos(v))^{3}\sin(v)-468\,\sin(v){v}^{2}-7742\,(\cos(v))^{3}{v}^{3}+96\,\sin(v){v}^{4}\\
-360\,\cos(v)\sin(v)+24240\,(\cos(v))^{4}v-29032\,(\cos(v))^{5}{v}^{3}+27360\,\sin(v)(\cos(v))^{6}\\
+53008\,(\cos(v))^{7}{v}^{3}-3840\,\sin(v)(\cos(v))^{9}-6016\,(\cos(v))^{10}{v}^{3}\\
-13728\,(\cos(v))^{8}{v}^{3}+17280\,(\cos(v))^{8}v-6720\,\sin(v)(\cos(v))^{7}\\
+32832\,\sin(v)(\cos(v))^{8}{v}^{2}+23040\,(\cos(v))^{9}v-6400\,(\cos(v))^{11}{v}^{3}\\
+13440\,\sin(v)(\cos(v))^{10}{v}^{2}-4608\,\sin(v)(\cos(v))^{7}{v}^{4}-4608\,\sin(v)(\cos(v))^{8}{v}^{4}\\
-11520\,\sin(v)(\cos(v))^{8}-7680\,\sin(v)(\cos(v))^{10}+10752\,(\cos(v))^{10}v\\
-1536\,\sin(v)(\cos(v))^{10}{v}^{4}+15360\,(\cos(v))^{11}v
\end{array}
\end{equation*}

Since for small values of $v$, the above formulae are subject to
heavy cancelations, we give the Taylor expansions of the
coefficients $b$:

\ason
\begin{equation*}
\begin{array}{l}
\ds b_1 = {\frac {399187}{241920}}-{\frac
{262795}{912384}}\,{v}^{2}+{\frac {
17265277}{4358914560}}\,{v}^{4}-{\frac
{38566679}{38041436160}}\,{v}^{ 6}-{\frac
{52935007231}{426824913715200}}\,{v}^{8}
 \\
\ds b_2= -{\frac {17327}{8640}}+{\frac {262795}{114048}}\,{v}^{2}-{\frac {
2649128441}{4358914560}}\,{v}^{4}+{\frac {2650726483}{52306974720}}\,{
v}^{6}-{\frac {374485131133}{106706228428800}}\,{v}^{8}
 \\
\ds b_3= {\frac {597859}{60480}}-{\frac {1839565}{228096}}\,{v}^{2}+{\frac {
1110676079}{311351040}}\,{v}^{4}-{\frac {6919361527}{8047226880}}\,{v}
^{6}+{\frac {1637603830619}{15243746918400}}\,{v}^{8}
 \\
\ds b_4= -{\frac {704183}{60480}}+{\frac {1839565}{114048}}\,{v}^{2}-{\frac {
5518849841}{622702080}}\,{v}^{4}+{\frac {14263211663}{4755179520}}\,{v
}^{6}-{\frac {73005517242211}{106706228428800}}\,{v}^{8}
\\
\ds b_5= {\frac {465133}{24192}}-{\frac {9197825}{456192}}\,{v}^{2}+{\frac {
734695627}{62270208}}\,{v}^{4}-{\frac {183227481067}{41845579776}}\,{v
}^{6}+{\frac {9908801489731}{8536498274304}}\,{v}^{8}
\end{array}
\end{equation*}
\asoff

The PLTE of the method is given

\begin{equation}
\begin{array}{l}
\ds plte=( {\frac {52559}{912384}}\,{w}^{10}y^{(2)}+{\frac {262795}{
456192}}\,{w}^{6}y^{(6)}+\nonumber \\
\ds {\frac{262795}{912384}}\,{w}^{8}y^{(4)}+{ \frac
{262795}{456192}}\,{w}^{4}y^{(8)}+{\frac {52559}{912384}}\,y^{(
12)}+{\frac {262795}{912384}}\,{w}^{2}y^{(10)} ) {h}^{12}
\end{array}
\end{equation}

\section{Stability Analysis}
The stability of the new methods is studied by considering the test equation
\begin{equation}
 \frac{d^2y(t)}{dt^2}=-\sigma ^2 y(t)
\end{equation}
and the linear multistep method (\ref{equ_ref_meth}) for the numerical solution.
In the above equation $\sigma \neq \omega$ ($\omega$ is the frequency at which the phase-lag function and its derivatives vanish). Setting $s=\sigma h$ and $v=\omega h$ we get for the characteristic equation of the applied method
\begin{equation}
 \sum_{j=1}^{J/2} A_j(s^2,v)(z^j+z^{-j})+A_0(s^2,v)=0
\end{equation}
where
\begin{equation}
 A_j(s^2,v)=a_{\frac{J}{2}-j}(v)+s^2\cdot b_{\frac{J}{2}-j}(v)
\end{equation}
The motivation of the above analysis is straightforward: Although the coefficients of the method (\ref{equ_ref_meth}) are designed in a way that the phase-lag and its first derivatives vanish in the frequency $\omega$, the frequency $\omega$ itself is unknown and only an estimation can be made. Thus, if the correct frequency of the problem is $\sigma$ we have to check if the method is stable, that is if the roots of the characteristic equation lie on the unit disk. For this reason we draw at the $v$-$s$ plane the areas in which the method is stable. Figure \ref{fig_res_sr} shows the stability region for the six methods (the classical one, the phase fitted one and those with first, second, third and fourth phase lag derivative elimination). Note here that the $s$-axis corresponds to the real frequency while the $v$-axis corresponds to the estimated frequency used to construct the parameters of the method.

\section{Numerical Results}
Numerical experiments have been carried out for two orbital problems. Since the classical method is well studied, we only present the new methods in comparison to the classical one.

\subsection{The 2-Body Problem}
In this problem we test the motion of two bodies in a reference system that is fixed in one of them. Moreover, the motion is planar, thus, we only have to calculate the $x$ and $y$ coordinates of the second body. The differential equations are:
\begin{eqnarray}
 \ddot{x}=&-\frac{x}{(x^2+y^2)^3} \nonumber \\
 \ddot{y}=&-\frac{y}{(x^2+y^2)^3} \nonumber
\end{eqnarray}
and the initial conditions are
\begin{eqnarray}
 x(0)=1-\epsilon & \dot{x}(0)=0 \nonumber \\
y(0)=0 & \dot{y}(0)=\sqrt{\frac{1+\epsilon}{1-\epsilon}} \nonumber
\end{eqnarray}
where $\epsilon$ is the eccentricity. Figure \ref{fig_res_2b} presents the accuracy of the methods expressed by $- log_{10}$ (error at
the end point which is after $100$ periods) versus the $log_{10}$ (total steps). The results clearly show the improvement in the accuracy using the proposed methodology.

\subsection{The 5-Outer Planet System}
The next problem concerns the motion of the five outer planets relative to the sun. The problem falls in the category of the N-Body problem which is the problem that regards the movement of N bodies under Newton's law of gravity. It is expressed by a system of vector differential
equations
\begin{equation}
\ddot{\bold{y}}_i=G\sum_{j=1,j\neq i}^N \frac{m_j(\bold{y_j}-\bold{y_i})}{|\bold{y_j}-\bold{y_i}|^3}
\label{equ_N_body}\end{equation}
where $G$ is the gravitational constant, $m_j$ is the mass of body $j$ and $\bold{y_i}$ is the vector of the position of body $i$. It is easy to see that each vector differential equation of (\ref{equ_N_body}) can be analyzed
into three simplified differential equations, that express the three directions
$x,y,z$.

The above system of ODEs cannot be solved analytically. Instead we produce a highly accurate numerical solution by using a 10-stage implicit Runge-Kutta method of Gauss with 20th algebraic order, that is also symplectic and A-stable. The method can be easily reproduced using simplifying assumptions for the order conditions (see \cite{butcher_Book_NMODE_W_03}). The reference solution is obtained by using the previous method to integrate the N-body problem for a specific time-span and for different step-lengths. In \cite{hairer_Book_GNISPAODE_S_02} the data for the five outer planet problem is given (these data are summarized in table \ref{table_5op}). Masses are relative to the sun, so that the sun has mass 1. In the computations the sun with the four inner planets are considered one body, so the mass is larger than one. Distances are in astronomical units, time is in earth days and the gravitational constant is $G = 2.95912208286 \cdot 10^{-4}$. The system of equations (\ref{equ_N_body}) has been solved for $t \in [0, 10^6 ]$, for which time-span, the previously mentioned method of Gauss produces a 10.5 decimal digits solution. Figure \ref{fig_res_5o} presents the accuracy of the methods expressed by $- log_{10}$ (error at the end point) versus the $log_{10}$ (total steps). The results clearly show the improvement in the accuracy using the proposed methodology.

\subsection{The figures}
\begin{figure}[hbt]
\includegraphics[width=\textwidth]{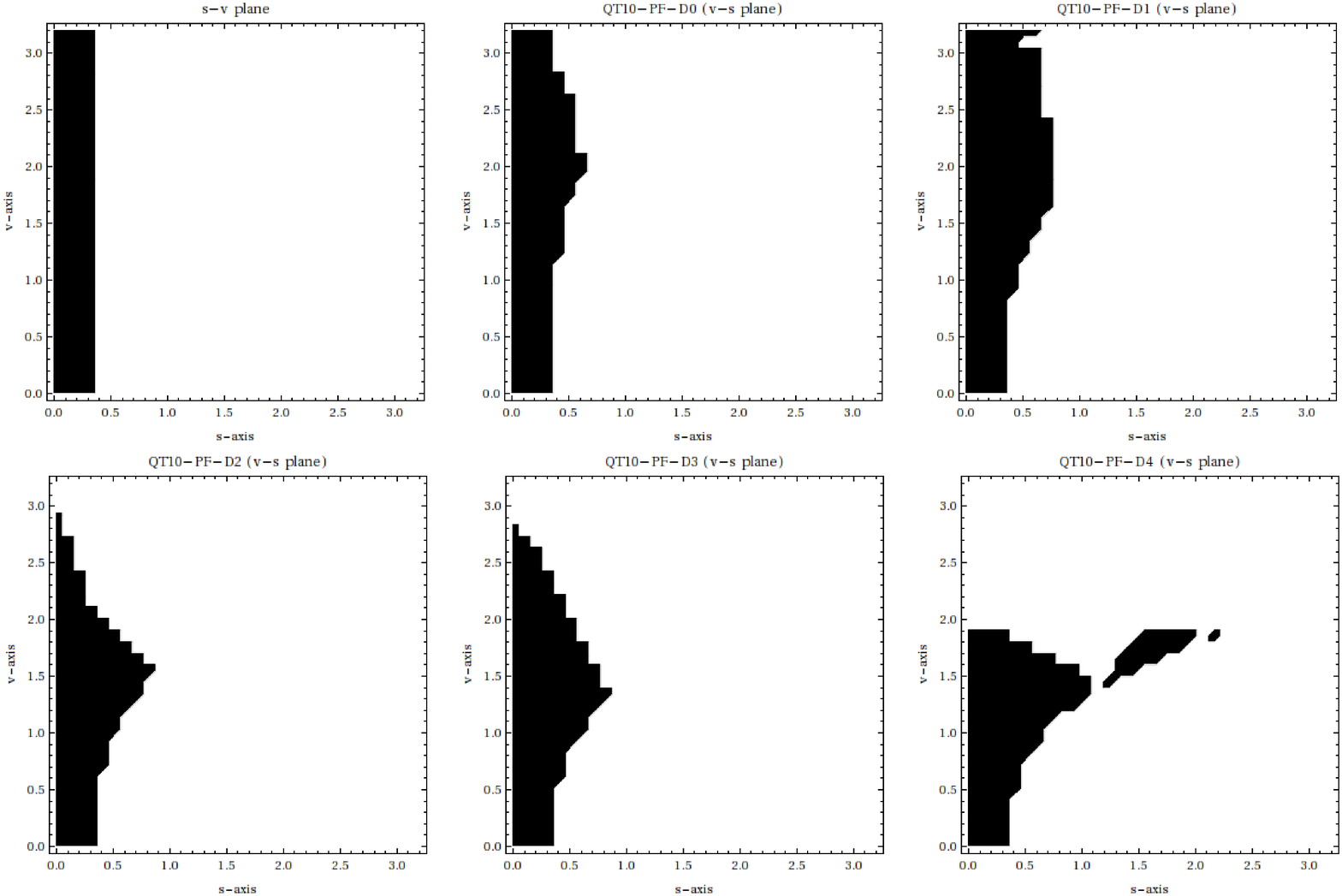}
\caption{Figure 1.}\label{fig_res_sr}
\end{figure}

\begin{figure}[hbt]
\includegraphics[width=\textwidth]{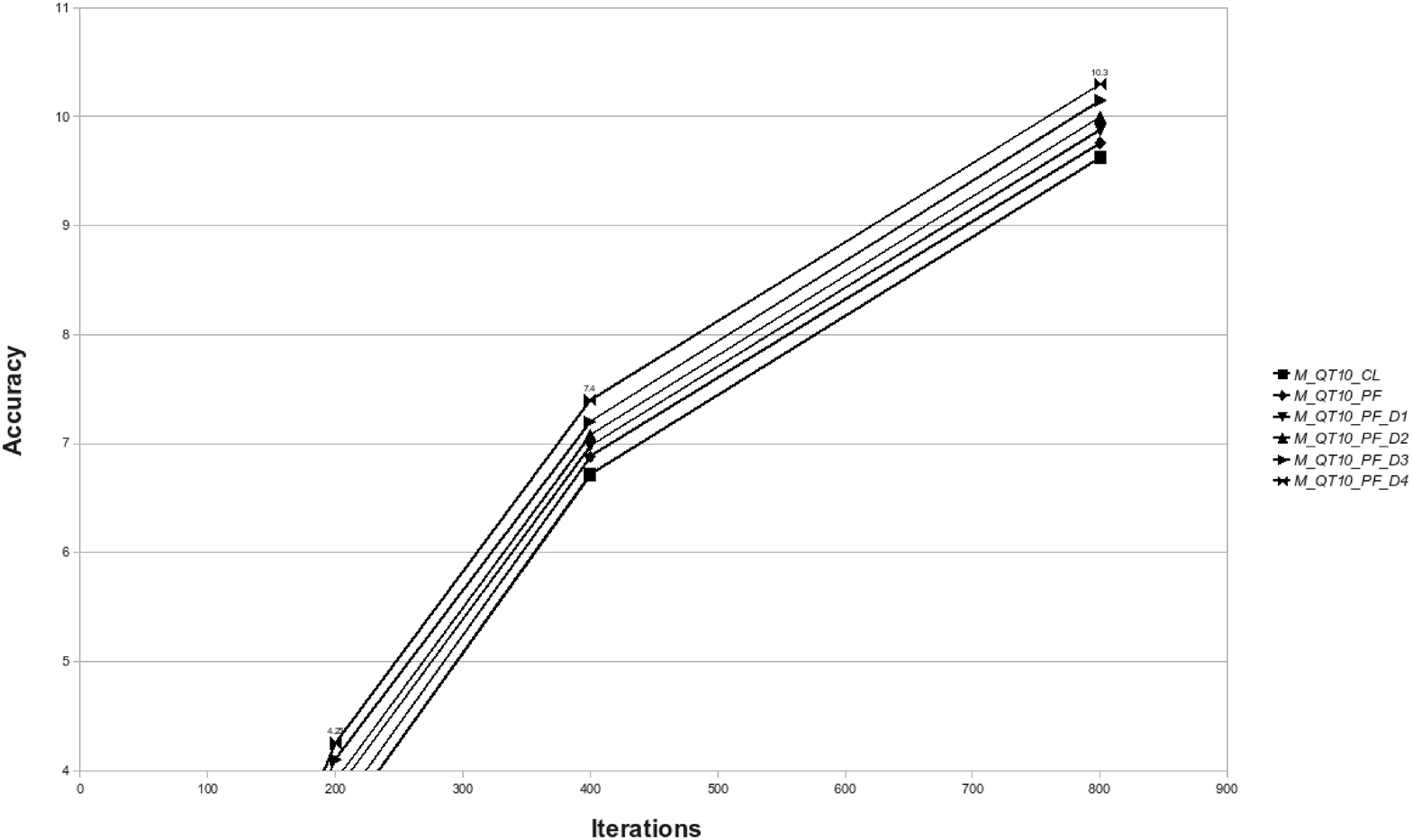}
\caption{Figure 2.}\label{fig_res_2b}
\end{figure}

\begin{figure}[hbt]
\includegraphics[width=\textwidth]{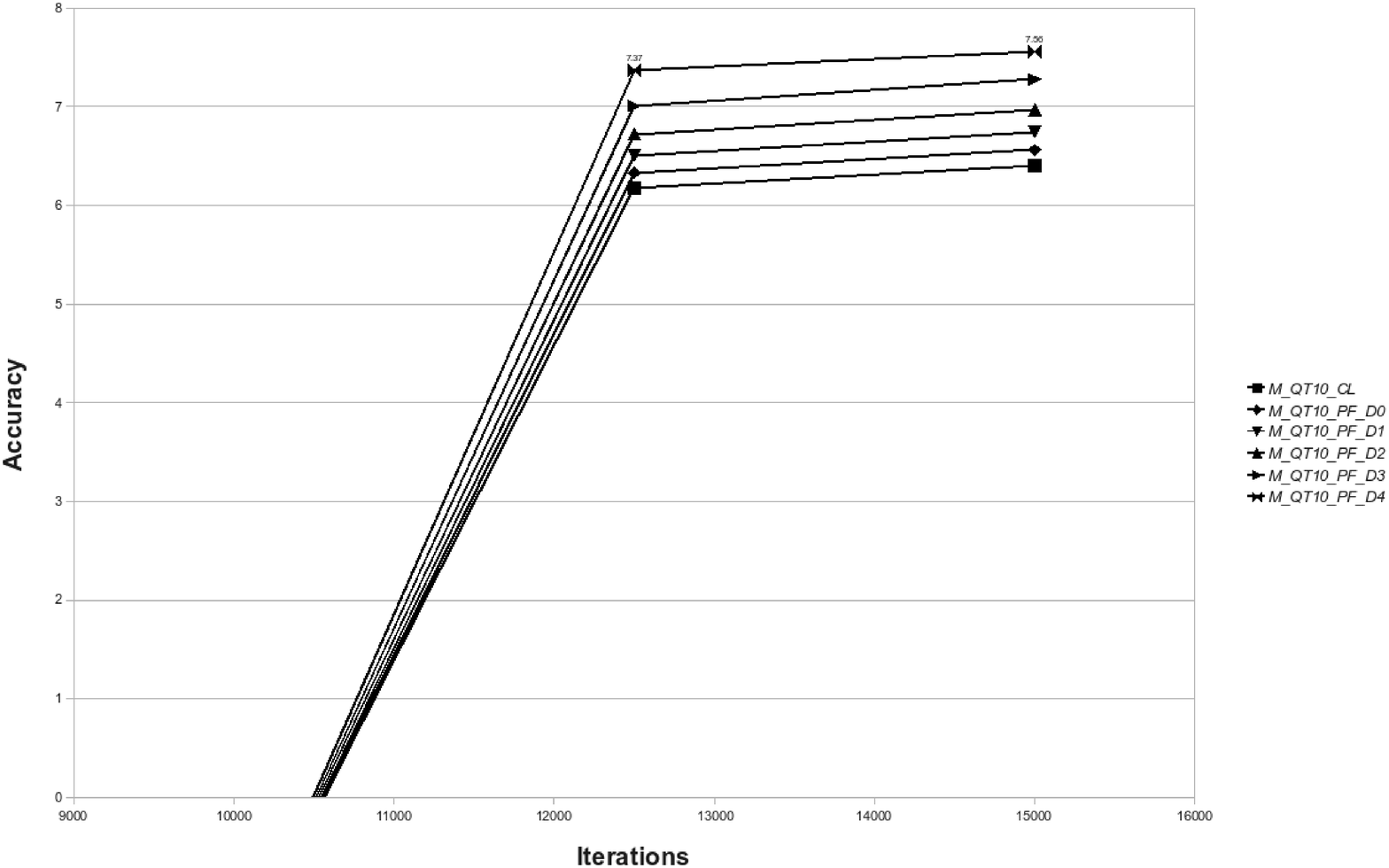}
\caption{Figure 3.}\label{fig_res_5o}
\end{figure}

In Figure \ref{fig_res_sr} we see the stability regions ($v$-$s$ plane) for the classical Quinlan-Tremaine method, the PF-D0, PF-D1, PF-D2, PF-D3 and PF-D4 methods (from left to right and from top to bottom). In Figure \ref{fig_res_2b} we see the accuracy of the methods expressed by $- log_{10}$ (error at the end point which is after $100$ periods) versus the $log_{10}$ (total steps) in the two body problem for eccentricity $\epsilon=0.5$ and step size $h=0.1$. In Figure \ref{fig_res_5o} we see the accuracy of the methods expressed by $- log_{10}$ (error at the end point which is after $10^6$ days) versus the $log_{10}$ (total steps) for the 5-Outer Planet System.

\begin{table}[htb]
\caption{Initial Data for the 5-outer planet problem}
\centering
\begin{tabular}{c | c | c | c}
\hline
\hline
Planet & Mass & Initial Position & Initial Velocity \\
\hline
Sun & 1.00000597682 & 0 & 0\\
 & & 0 & 0\\
 & & 0 & 0\\
\hline
Jupiter & 0.000954786104043 & -3.5023653 & 0.00565429\\
 & & -3.8169847 & -0.00412490\\
 & & -1.5507963 & -0.00190589\\
\hline
Saturn & 0.000285583733151 & 9.0755314 & 0.00168318\\
 & & -3.0458353 & 0.00483525\\
 & & -1.6483708 & 0.00192462\\
\hline
Uranus & 0.0000437273164546 & 8.3101420 & 0.00354178\\
 & & -16.2901086 & 0.00137102\\
 & & -7.2521278 & 0.00055029\\
\hline
Neptune & 0.0000517759138449 & 11.4707666 & 0.00288930\\
 & & -25.7294829 & 0.00114527\\
 & & -10.8169456 & 0.00039677\\
\hline
Pluto & $1/(1.3 \cdot 10^8 )$ & -15.5387357 & 0.00276725\\
 & & -25.2225594 & -0.00170702\\
 & & -3.1902382 & -0.00136504
\end{tabular}
\label{table_5op}
\end{table}

\section{Conclusions}
We have presented a new family of 10-step symmetric multistep numerical methods with improved characteristics concerning orbital problems. The methods were constructed by adopting a new methodology which, apart from the phase fitting at a predefined frequency, it eliminates the first derivatives of the phase lag function at the same frequency. The result is that the phase lag function becomes less sensitive on the frequency near the predefined one. This behavior compensates the fact that the exact frequency can only be estimated. Experimental results demonstrate this behavior by showing that the accuracy is increased as the number of the derivatives that are eliminated is increased.

\bibliographystyle{apj}
\bibliography{../../../math}

\begin{thebibliography}{24}
\expandafter\ifx\csname natexlab\endcsname\relax\def\natexlab#1{#1}\fi

\bibitem[{Anastassi \& Simos(2004)}]{anastassi_IJMPC_15_1_04}
Anastassi, Z. \& Simos, T. 2004, International Journal of Modern Physics C, 15,
  1

\bibitem[{Anastassi \& Simos(2005{\natexlab{a}})}]{anastassi_MCM_42_877_05}
---. 2005{\natexlab{a}}, Mathematical and Computer Modelling, 42, 877

\bibitem[{Anastassi \& Simos(2005{\natexlab{b}})}]{anastassi_NA_10_301_05}
---. 2005{\natexlab{b}}, New Astronomy, 10, 301

\bibitem[{Anastassi \& Simos(2007)}]{anastassi_JMC_41_79_07}
---. 2007, Journal of Mathematical Chemistry, 41, 79

\bibitem[{Berghe \& Daele(2006)}]{berghe_JNAIAM_1_241_06}
Berghe, G.~V. \& Daele, M.~V. 2006, JNAIAM, 1, 241

\bibitem[{Brusa \& Nigro(1980)}]{brusa_IJNME_15_685_80}
Brusa, L. \& Nigro, L. 1980, Int. J. Num. Methods Engrg., 15, 685–699

\bibitem[{Butcher(2003)}]{butcher_Book_NMODE_W_03}
Butcher, J. 2003, Numerical methods for ordinary differential equations (Wiley)

\bibitem[{Cash \& Mazzia(2006)}]{cash_JNAIAM_1_81_06}
Cash, J.~R. \& Mazzia, F. 2006, JNAIAM, 1, 89

\bibitem[{Chawla \& Rao(1986)}]{chawla_JCAM_15_329_86}
Chawla, M. \& Rao, P. 1986, J.Comput.Appl.Math., 15, 329

\bibitem[{Gautschi(1961)}]{gautschi_NM_3_381_61}
Gautschi, W. 1961, Numer. Math., 3, 381–397

\bibitem[{Hairer {et~al.}(2002)Hairer, Lubich, \&
  Wanner}]{hairer_Book_GNISPAODE_S_02}
Hairer, E., Lubich, C., \& Wanner, G. 2002, Geometric numerical integration,
  Structure preserving algorithms for ordinary differential equations
  (Springer)

\bibitem[{Iavernaro {et~al.}(2006)Iavernaro, Mazzia, \&
  Trigiante}]{iavernaro_JNAIAM_1_91_06}
Iavernaro, F., Mazzia, F., \& Trigiante, D. 2006, JNAIAM, 1, 91

\bibitem[{Ixaru \& Berghe(2004)}]{ixaru_Book_EF_KAP_04}
Ixaru, L. \& Berghe, G.~V. 2004, Exponential Fitting (Dordrecht/Boston/London:
  Kluwer Academic Publishers)

\bibitem[{Lambert \& Watson(1976)}]{lambert_JIMA_18_189_76}
Lambert, J. \& Watson, I. 1976, J. Inst. Math. Appl., 18, 189

\bibitem[{Lyche(1972)}]{lyche_NM_19_65_72}
Lyche, T. 1972, Num. Math., 19, 65

\bibitem[{Mazzia {et~al.}(2006)Mazzia, Sestini, \&
  Trigiante}]{mazzia_JNAIAM_1_131_06}
Mazzia, F., Sestini, A., \& Trigiante, D. 2006, JNAIAM, 1, 131

\bibitem[{Psihoyios(2006)}]{psihoyios_CL_2_51_06}
Psihoyios, G. 2006, CoLe, 2, 51

\bibitem[{Quinlan \& Tremaine(1990)}]{quinlan_AJ_100_1694_90}
Quinlan, D. \& Tremaine, S. 1990, The Astronomical Journal, 100, 1694

\bibitem[{Quinlan(1999)}]{quinlan_arxiv_astro_ph_9901136}
Quinlan, G. 1999, preprint arXiv, astro-ph/9901136

\bibitem[{Raptis \& Allison(1978)}]{raptis_CPC_14_1_78}
Raptis, D. \& Allison, A. 1978, Computer Physics Communications, 14, 1

\bibitem[{Simos(2000)}]{simos_Book_CMATV1_RCS_00}
Simos, T. 2000, Specialist Periodical Reports (Cambridge: The Royal Society of
  Chemistry)

\bibitem[{Simos(2005)}]{simos_CL_1_37_05}
---. 2005, CoLe, 1, 37

\bibitem[{Simos(2007)}]{simos_CL_1_45_07}
---. 2007, CoLe, 1, 45

\bibitem[{Simos \& Williams(1999)}]{simos_CC_23_513_99}
Simos, T. \& Williams, P. 1999, Comput. Chem., 23, 513

\end{thebibliography}

\end{document}